\documentclass[11pt]{amsart}
\usepackage{amsmath,amssymb,mathrsfs}
\newtheorem{theorem}{Theorem}

\newtheorem{lemma}[theorem]{Lemma}

\theoremstyle{definition}

\begin{document}

\title{Alexandrov immersed minimal tori in $S^3$}
\author{Simon Brendle}
\begin{abstract}
We show that any minimal torus in $S^3$ which is Alexandrov immersed must be rotationally symmetric. An analogous result holds for surfaces of constant mean curvature.
\end{abstract}
\address{Department of Mathematics \\ Stanford University \\ Stanford, CA 94305}
\thanks{The author was supported in part by the National Science Foundation under grants DMS-0905628 and DMS-1201924.}
\maketitle

\section{Introduction}

In a recent paper \cite{Brendle1}, we showed the Clifford torus is the only embedded minimal surface in $S^3$ of genus $1$, thereby confirming a conjecture of Lawson. In this note, we classify minimal tori in $S^3$ that are immersed in the sense of Alexandrov.

\begin{theorem}
\label{alex.immersed}
Let $F: \Sigma \to S^3$ be an immersed minimal surface in $S^3$ of genus $1$. Moreover, we assume that $F$ is an Alexandrov immersion; this means that there exists a compact manifold $N$ and an immersion $\bar{F}: N \to S^3$ such that $\partial N = \Sigma$ and $\bar{F}|_\Sigma = F$. Then $\Sigma$ is rotationally symmetric. That us, we can find an anti-symmetric matrix $Q \in \mathfrak{so}(4)$ of rank $2$ such that $Q \, F(x) \in \text{\rm span} \{\frac{\partial F}{\partial x_1}(x),\frac{\partial F}{\partial x_2}(x)\}$ for all $x \in \Sigma$.
\end{theorem} 

There is a complete classification of all rotationally symmetric minimal tori in $S^3$; for details, we refer to \cite{Perdomo} or \cite{Brendle3}, Theorem 1.4. Besides the Clifford torus, there is a large class of additional examples which are Alexandrov immersed but fail to be embedded. 

We will present the proof of Theorem \ref{alex.immersed} in Section \ref{minimal}. The argument is similar in spirit to the case of embedded surfaces studied in \cite{Brendle1}, and we will only indicate the necessary modifications.

After the paper \cite{Brendle1} was published, Andrews and Li \cite{Andrews-Li} observed that the arguments in \cite{Brendle1} can be extended to the setting of constant mean curvature surfaces. As a result, they showed that every embedded constant mean curvature surface in $S^3$ is rotationally symmetric. More recently, we showed that the arguments in \cite{Brendle1} can be generalized to a class of Weingarten surfaces (see \cite{Brendle2}).

Our proof of Theorem \ref{alex.immersed} also extends to the setting of constant mean curvature surfaces. This yields the following result:

\begin{theorem} 
\label{cmc}
Let $F: \Sigma \to S^3$ be an immersed constant mean curvature surface in $S^3$ of genus $1$. Suppose that $F$ extends to an immersion $\bar{F}: N \to S^3$ where $\partial N = \Sigma$ and that $\partial N$ is mean convex with respect to the pull-back of the standard metric on $S^3$ under $\bar{F}$. Then $\Sigma$ is rotationally symmetric.
\end{theorem}

The proof of Theorem \ref{cmc} is similar to Theorem \ref{alex.immersed}. The condition that the surface is Alexandrov immersed is quite natural in light of the work of Korevaar, Kusner, and Ratzkin \cite{Korevaar-Kusner-Ratzkin}, Korevaar, Kusner, and Solomon \cite{Korevaar-Kusner-Solomon}, and Kusner, Mazzeo, and Pollack \cite{Kusner-Mazzeo-Pollack}, where Alexandrov immersed constant mean curvature surfaces in $\mathbb{R}^3$ have been studied.

The author is very grateful to Professor Robert Kusner for bringing this problem to his attention, and for many helpful discussions. He would like to thank the Department of Mathematics at Cambridge University, where part of this work was carried out.

\section{Proof of Theorem \ref{alex.immersed}}

\label{minimal}

For convenience, we put a Riemannian metric on $N$ so that $\bar{F}$ is a local isometry. In particular, there exists a real number $\delta>0$ so that $\bar{F}(x) \neq \bar{F}(y)$ for all points $x,y \in N$ satisfying $d_N(x,y) \in (0,\delta)$.

For each point $x \in \Sigma$ and any number $\alpha \geq 1$, we define 
\[D_\alpha(x) = \Big \{ p \in S^3: \frac{\alpha}{\sqrt{2}} \, |A(x)| \, (1 - \langle F(x),p \rangle) + \langle \nu(x),p \rangle \leq 0 \Big \}.\] 
Note that $D_\alpha(x)$ is a geodesic ball in $S^3$ whose boundary passes through the point $F(x)$. Moreover, the outward-pointing unit normal vector to $\partial D_\alpha(x)$ at the point $F(x)$ is given by $\nu(x)$. 

Let $I$ denote the set of all points $(x,\alpha) \in \Sigma \times [1,\infty)$ with the property that there exists a smooth map $G: D_\alpha(x) \to N$ such that $\bar{F} \circ G = \text{\rm id}_{D_\alpha(x)}$ and $G(F(x)) = x$. 

\begin{lemma}
\label{lemma.1}
Let us fix a pair $(x,\alpha) \in I$. Then there is a unique map $G: D_\alpha(x) \to N$ such that $\bar{F} \circ G = \text{\rm id}_{D_\alpha(x)}$ and $G(F(x)) = x$.
\end{lemma}

\textbf{Proof.} 
Suppose that there are two maps $G,\tilde{G}: D_\alpha(x) \to N$ such that $\bar{F} \circ G = \bar{F} \circ \tilde{G} = \text{\rm id}_{D_\alpha(x)}$ and $G(F(x)) = \tilde{G}(F(x)) = x$. For each point $p \in D_\alpha(x)$, we have $\bar{F}(G(p)) = \bar{F}(\tilde{G}(p))$, hence $d_N(G(p),\tilde{G}(p)) \notin (0,\delta)$. By continuity, we either have $G(p) = \tilde{G}(p)$ for all $p \in D_\alpha(x)$ or we have $G(p) \neq \tilde{G}(p)$ for all $p \in D_\alpha(x)$. The second case can be ruled out, as $G(F(x)) = \tilde{G}(F(x))$. Thus, we conclude that $G(p) = \tilde{G}(p)$ for all $p \in D_\alpha(x)$. This shows that $G$ is unique. \\

\begin{lemma}
\label{lemma.2}
The set $I$ is closed. Furthermore, the map $G$ depends continuously on the pair $(x,\alpha)$.
\end{lemma}

\textbf{Proof.} 
Consider a sequence of pairs $(x^{(m)},\alpha^{(m)}) \in I$ such that $\lim_{m \to \infty} (x^{(m)},\alpha^{(m)}) = (\hat{x},\hat{\alpha})$. By definition of $I$, we can find a sequence of smooth maps $G^{(m)}: D_{\alpha^{(m)}}(x^{(m)}) \to N$ such that $\bar{F} \circ G^{(m)} = \text{\rm id}_{D_{\alpha^{(m)}}(x^{(m)})}$ and $G^{(m)}(F(x^{(m)})) = x^{(m)}$. Since $\bar{F}$ is a smooth immersion, the maps $G^{(m)}$ are uniformly bounded in $C^2$ norm. Hence, after passing to a subsequence, the maps $G^{(m)}$ converge in $C^1$ to a map $G: D_{\hat{\alpha}}(\hat{x}) \to N$. The map $G$ satisfies $\bar{F} \circ G = \text{\rm id}_{D_{\hat{\alpha}}(\hat{x})}$ and $G(F(\hat{x})) = \hat{x}$. From this, we deduce that $G$ is smooth. Thus, $(\hat{x},\hat{\alpha}) \in I$, and the assertion follows. \\

\begin{lemma}
\label{lemma.3}
We have $(x,\alpha) \in I$ if $\alpha$ is sufficiently large. 
\end{lemma}

\textbf{Proof.} 
By a result of Lawson, we have $|A(x)| > 0$ for all $x \in \Sigma$. Hence, if we choose $\alpha$ sufficiently large, the radius of the geodesic ball $D_\alpha(x) \subset S^3$ can be made arbitrarily small. Therefore, if $\alpha$ is sufficiently large, the implicit function theorem guarantees the existence of a smooth map $G: D_\alpha(x) \to N$ such that $\bar{F} \circ G = \text{\rm id}_{D_\alpha(x)}$ and $G(F(x)) = x$. \\

After these preparations, we now describe the proof of Theorem \ref{alex.immersed}. Let us define 
\[\kappa = \inf \{\alpha: \text{\rm $(x,\alpha) \in I$ for all $x \in \Sigma$}\}.\] 
Clearly, $\kappa \in [1,\infty)$. For each point $x \in \Sigma$, there is a unique map $G_x: D_\kappa(x) \to N$ such that $\bar{F} \circ G_x = \text{\rm id}_{D_\kappa(x)}$ and $G_x(F(x)) = x$. For each point $x \in \Sigma$, the map $G_x$ and the map $\bar{F}|_{G_x(D_\kappa(x))}$ are one-to-one. 

We next define a smooth function $Z: \Sigma \times \Sigma \to \mathbb{R}$ by 
\[Z(x,y) = \frac{\kappa}{\sqrt{2}} \, |A(x)| \, (1 - \langle F(x),F(y) \rangle) + \langle \nu(x),F(y) \rangle\] 
for $x,y \in \Sigma$. In contrast to \cite{Brendle1}, the function $Z(x,y)$ might be negative somewhere.

As in \cite{Brendle1}, we distinguish two cases: \\

\textbf{Case 1:} Suppose first that $\kappa = 1$. 

\begin{lemma} 
We have $Z(x,y) \geq 0$ if $x$ and $y$ are sufficiently close.
\end{lemma}

\textbf{Proof.} 
The proof is by contradiction. Suppose that there exist two sequences of points $x^{(m)},y^{(m)} \in \Sigma$ such that $\lim_{m \to \infty} x^{(m)} = \lim_{m \to \infty} y^{(m)}$ and $Z(x^{(m)},y^{(m)}) < 0$ for all $m$. Since $Z(x^{(m)},y^{(m)}) < 0$, the point $F(y^{(m)})$ lies in the interior of the geodesic ball $D_\kappa(x^{(m)})$. Since $G_{x^{(m)}}$ is an immersion, the point $\tilde{y}^{(m)} := G_{x^{(m)}}(F(y^{(m)}))$ must lie in the interior of $N$. Since $y^{(m)}$ lies on the boundary $\Sigma$, it follows that 
\[\tilde{y}^{(m)} \neq y^{(m)}.\] 
On the other hand, we have 
\[\bar{F}(\tilde{y}^{(m)}) = F(y^{(m)})\] 
and 
\[\lim_{m \to \infty} \tilde{y}^{(m)} = \lim_{m \to \infty} G_{x^{(m)}}(F(y^{(m)})) = \lim_{m \to \infty} G_{x^{(m)}}(F(x^{(m)})) = \lim_{m \to \infty} x^{(m)} = \lim_{m \to \infty} y^{(m)}.\] 
This is impossible since $\bar{F}$ is an immersion. \\

\begin{lemma}
\label{gradient.of.curvature}
Fix a point $x \in \Sigma$, and let $\{e_1,e_2\}$ be an orthonormal basis of $T_x \Sigma$ such that $h(e_1,e_1) > 0$, $h(e_1,e_2) = 0$, and $h(e_2,e_2) < 0$. Then $e_1(|A|) = 0$.
\end{lemma}

\textbf{Proof.} 
The proof is analogous to arguments in \cite{Brendle1}. Let $\gamma: \mathbb{R} \to \Sigma$ be a geodesic such that $\gamma(0) = x$ and $\gamma'(0) = e_1$. Since the function $Z$ is nonnegative when $x$ and $y$ are sufficiently close, the function $f(t) = Z(x,\gamma(t))$ is nonnegative when $t$ is sufficiently small. Since $\kappa=1$, we have $f(0) = f'(0) = f''(0) = 0$. Since the function $f(t)$ is nonnegative in a neighborhood of $0$, we conclude that $f'''(0) = 0$. From this, we deduce that $(D_{e_1}^\Sigma h)(e_1,e_1) = 0$. From this, the assertion follows. \\

Lemma \ref{gradient.of.curvature} implies that the function $|A|$ is constant along one family of curvature lines on $\Sigma$. We claim that $F$ is rotationally symmetric. To see this, we define a vector field $V$ on $\Sigma$ by $V = |A|^{-\frac{1}{2}} \, e_1$. The identity $e_1(|A|) = 0$ implies that $[V,e_1] = 0$. Moreover, we compute 
\[D_{e_1}^\Sigma e_2 = -\frac{1}{2 \, |A|} \, e_2(|A|) \, e_1\] 
and 
\[D_{e_2}^\Sigma e_1 = 0.\] 
This implies 
\[[e_1,e_2] = D_{e_1}^\Sigma e_2 - D_{e_2}^\Sigma e_1 = -\frac{1}{2 \, |A|} \, e_2(|A|) \, e_1,\] 
hence 
\[[V,e_2] = |A|^{-\frac{1}{2}} \, [e_1,e_2] - e_2(|A|^{-\frac{1}{2}}) \, e_1 = 0.\] 
Using the identities $[V,e_1] = [V,e_2] = 0$, it is straightforward to check that $\mathscr{L}_V g = \mathscr{L}_V h = 0$. Therefore, $V$ is the restriction of an ambient Killing vector field on $S^3$. Hence, we can find a constant matrix $Q \in \mathfrak{so}(4)$ such that $V(x) = Q \, F(x)$ for all $x \in \Sigma$. We now differentiate this relation along $e_2$. Since $g(e_2,V) = h(e_2,V) = 0$, we obtain the relation $D_{e_2}^\Sigma V = Q \, e_2$. Since $D_{e_2}^\Sigma V = e_2(|A|^{-\frac{1}{2}}) \, e_1$, we conclude that 
\[Q \, \Big ( e_2 + \frac{1}{2 \, |A|} \, e_2(|A|) \, F(x) \Big ) = Q \, e_2 + \frac{1}{2 \, |A|^{\frac{3}{2}}} \, e_2(|A|) \, e_1 = 0.\] 
Therefore, the matrix $Q$ has non-trivial nullspace. Thus, $Q$ has rank $2$ and $F$ is rotationally symmetric. \\

\textbf{Case 2:} Suppose next that $\kappa > 1$. 

\begin{lemma} 
\label{touch.1}
Fix a point $x \in \Sigma$. Then there exists a constant $\beta > 0$ such that $d_N(G_x(p),\Sigma) \geq \beta \, |p - F(x)|^2$ for all points $p \in \partial D_\kappa(x)$ that are sufficiently close to $F(x)$.
\end{lemma}

\textbf{Proof.} 
Let us fix a point $x \in \Sigma$. Moreover, we consider the function 
\[\varphi_x: \partial D_\kappa(x) \to \mathbb{R}, \quad p \mapsto d_N(G_x(p),\Sigma).\] 
Clearly, $\varphi_x(F(x)) = 0$, and the gradient of the function $\varphi_x$ at the point $F(x)$ vanishes. Moreover, since $\kappa > 1$, the Hessian of the function $\varphi_x$ at the point $F(x)$ is positive definite. Hence, we can find a positive constant $\beta > 0$ such that $\varphi_x(p) \geq \beta \, |p - F(x)|^2$ for all points $p \in \partial D_\kappa(x)$ that are sufficiently close to $F(x)$. \\

\begin{lemma} 
\label{touch.2}
There exists a point $\hat{x} \in \Sigma$ such that $\Sigma \cap G_{\hat{x}}(D_\kappa(\hat{x})) \neq \{\hat{x}\}$.
\end{lemma}

\textbf{Proof.} 
Suppose this is false. Then $\Sigma \cap G_x(D_\kappa(x)) = \{x\}$ for all $x \in \Sigma$. This implies that $d_N(G_x(p),\Sigma) > 0$ for all $x \in \Sigma$ and all points $p \in \partial D_\kappa(x) \setminus \{F(x)\}$. Using Lemma \ref{touch.1}, we conclude that there exists a positive constant $\gamma > 0$ such that $d_N(G_x(p),\Sigma) \geq \gamma \, |p - F(x)|^2$ for all $x \in \Sigma$ and all $p \in \partial D_\kappa(x)$. By the implicit function theorem, there exists a small number $\varepsilon > 0$ such that $(x,\kappa-\varepsilon) \in I$ for all $x \in \Sigma$. This contradicts the definition of $\kappa$. \\

Let $\hat{x}$ be chosen as in Lemma \ref{touch.2}. Moreover, let us pick a point $\hat{y} \in \Sigma \cap G_{\hat{x}}(D_\kappa(\hat{x}))$ such that $\hat{x} \neq \hat{y}$. Since $\hat{y} \in G_{\hat{x}}(D_\kappa(\hat{x}))$, we conclude that $F(\hat{y}) \in D_\kappa(\hat{x})$ and $G_{\hat{x}}(F(\hat{y})) = \hat{y}$. Moreover, we claim that $F(\hat{x}) \neq F(\hat{y})$; indeed, if $F(\hat{x}) = F(\hat{y})$, then $\hat{x} = G_{\hat{x}}(F(\hat{x})) = G_{\hat{x}}(F(\hat{y})) = \hat{y}$, which contradicts our choice of $\hat{y}$.

We next consider the function $Z$ defined above. We claim that the function $Z$ is nonnegative in a neighborhood of the point $(\hat{x},\hat{y})$.

\begin{lemma}
We have $Z(x,y) \geq 0$ if $(x,y)$ is sufficiently close to $(\hat{x},\hat{y})$. 
\end{lemma}

\textbf{Proof.} 
We argue by contradiction. Suppose that there exist sequences of points $x^{(m)},y^{(m)} \in \Sigma$ such that $\lim_{m \to \infty} x^{(m)} = \hat{x}$, $\lim_{m \to \infty} y^{(m)} = \hat{y}$, and $Z(x^{(m)},y^{(m)}) < 0$. Since $Z(x^{(m)},y^{(m)}) < 0$, the point $F(y^{(m)})$ lies in the interior of the ball $D_\kappa(x^{(m)})$. Since $G_{x^{(m)}}$ is an immersion, we conclude that the point $\tilde{y}^{(m)} := G_{x^{(m)}}(F(y^{(m)}))$ lies in the interior of $N$. In particular, we have 
\[\tilde{y}^{(m)} \neq y^{(m)}.\] 
On the other hand, we have 
\[\bar{F}(\tilde{y}^{(m)}) = F(y^{(m)})\] 
and 
\[\lim_{m \to \infty} \tilde{y}^{(m)} = \lim_{m \to \infty} G_{x^{(m)}}(F(y^{(m)})) = G_{\hat{x}}(F(\hat{y})) = \hat{y} = \lim_{m \to \infty} y^{(m)}.\] 
This contradicts the fact that $\bar{F}$ is an immersion. Thus, $Z(x,y) \geq 0$ for $(x,y)$ close to $(\hat{x},\hat{y})$. \\

Therefore, we can find disjoint open sets $U,V \subset \Sigma$ such that $\hat{x} \in U$, $\hat{y} \in V$, and $Z(x,y) \geq 0$ for all points $(x,y) \in U \times V$. As in \cite{Brendle1}, we define 
\[\Omega = \{x \in U: \text{\rm there exists a point $y \in V$ such that $Z(x,y) = 0$}\}.\] 
Since $F(\hat{y}) \in D_\kappa(\hat{x})$, we have $Z(\hat{x},\hat{y}) \leq 0$. This implies that $Z(\hat{x},\hat{y}) = 0$, hence $\hat{x} \in \Omega$. We can now use the calculation in \cite{Brendle1} to conclude that $Z$ is a supersolution of a degenerate elliptic equation. More precisely, suppose that $(\bar{x},\bar{y})$ is an arbitrary point in $U \times V$. Then we can find a system of geodesic normal coordinates $(x_1,x_2)$ around $\bar{x}$ and a system of geodesic normal coordinates $(y_1,y_2)$ around $\bar{y}$ such that 
\begin{align*} 
&\sum_{i=1}^2 \frac{\partial^2 Z}{\partial x_i^2}(\bar{x},\bar{y}) + 2 \sum_{i=1}^2 \frac{\partial^2 Z}{\partial x_i \, \partial y_i}(\bar{x},\bar{y}) + \sum_{i=1}^2 \frac{\partial^2 Z}{\partial y_i^2}(\bar{x},\bar{y}) \\ 
&\leq -\frac{\kappa^2-1}{\sqrt{2} \, \kappa} \, \frac{|A(\bar{x})|}{1-\langle F(\bar{x}),F(\bar{y}) \rangle} \sum_{i=1}^2 \Big \langle \frac{\partial F}{\partial x_i}(\bar{x}),F(\bar{y}) \Big \rangle^2 \\ 
&+ \Lambda \, \bigg ( Z(\bar{x},\bar{y}) + \sum_{i=1}^2 \Big | \frac{\partial Z}{\partial x_i}(\bar{x},\bar{y}) \Big | + \sum_{i=1}^2 \Big | \frac{\partial Z}{\partial y_i}(\bar{x},\bar{y}) \Big | \bigg ), 
\end{align*} 
where $\Lambda$ is a positive constant. Using Bony's version of the strict maximum principle, we conclude that the set $\Omega$ contains an open neighborhood of $\hat{x}$. Moreover, the gradient of $|A|$ vanishes on the set $\Omega$. By analytic continuation, $|A|$ is a constant function on $\Sigma$. This implies that $F$ is congruent to the Clifford torus. This completes the proof of Theorem \ref{alex.immersed}. \\

\section{Sketch of the proof of Theorem \ref{cmc}}

Finally, let us sketch the proof of Theorem \ref{cmc}. Let $F: \Sigma \to S^3$ be an immersed constant mean curvature surface in $S^3$ of genus $1$. We assume that $F$ extends to an immersion $\bar{F}: N \to S^3$ where $\partial N = \Sigma$ and that $\partial N$ is mean convex with respect to the pull-back of the standard metric on $S^3$ under $\bar{F}$. Given a point $x \in \Sigma$ and a real number $\alpha \geq 1$, one defines 
\[D_\alpha(x) = \Big \{ p \in S^3: \Big ( \frac{H}{2} + \frac{\alpha}{\sqrt{2}} \, |\mathring{A}(x)| \Big ) \, (1 - \langle F(x),p \rangle) + \langle \nu(x),p \rangle \leq 0 \Big \}\] 
(cf. \cite{Andrews-Li}). Here, $H$ is the mean curvature (i.e. the sum of the principal curvatures) and $\mathring{A}$ is the trace-free part of the second fundamental form. As above, let $I$ denote the set of all points $(x,\alpha) \in \Sigma \times [1,\infty)$ with the property that there exists a smooth map $G: D_\alpha(x) \to N$ such that $\bar{F} \circ G = \text{\rm id}_{D_\alpha(x)}$ and $G(F(x)) = x$. It is well known that a constant mean curvature torus in $S^3$ has no umbilic points, so $|\mathring{A}(x)| > 0$ for all points $x \in I$. Hence, if $\alpha$ is sufficiently large, then the radius of the geodesic ball $D_\alpha(x)$ will be very small. From this, we deduce that $(x,\alpha) \in I$ if $\alpha$ is sufficiently large. We then define 
\[\kappa = \inf \{\alpha: \text{\rm $(x,\alpha) \in I$ for all $x \in \Sigma$}\}.\] 
If $\kappa=1$, we can argue as above to conclude that $e_1(|A|) = 0$, where $e_1$ is one of the eigenvectors of the second fundamental form. This implies that $F$ is rotationally symmetric. On the other hand, if $\kappa > 1$, we can combine the arguments above with the calculations in \cite{Andrews-Li} and \cite{Brendle1} to arrive at a contradiction.

\end{document}